\documentclass[11pt,a4paper]{amsart}

\usepackage{amssymb}
\usepackage{array,dcolumn}
\newcolumntype{C}{>{$}c<{$}}
\newcolumntype{L}{>{$}l<{$}}
\newcolumntype{R}{>{$}r<{$}}

\newtheorem*{theorem}{Theorem A}

\newcommand{\Q}{\mathbb{Q}}
\newcommand{\om}{\omega}

\renewcommand{\)}{{)\!)}}

\newcommand{\CP}{\mathbb{CP}}
\newcommand{\Mbar}{\overline{\mathcal{M}}}

\begin{document}

\title{Appendix to ``Intersection theory on $\Mbar_{1,4}$ and elliptic
Gromov-Witten invariants''}

\author{E. Getzler}

\address{Max-Planck-Institut f\"ur Mathematik, Gottfried-Claren-Str.\ 26,
D-53225 Bonn, Germany}

\curraddr{Department of Mathematics, Northwestern University, Evanston, IL
60208-2730, USA}

\email{getzler@math.nwu.edu}

\maketitle

\noindent
\textbf{Appendix: The recursions for $\CP^3$.}
\medskip

The genus $g$ Gromov-Witten potentials of the projective space $\CP^3$ have
the form
$$
F_g(\CP^3) = \begin{cases}
\bigl( \frac{1}{2} t_0^2t_3 + t_0t_1t_2 + \frac{1}{6} t_1^3 \bigr) +
\displaystyle \sum_{4\beta=a+2b} N^{(0)}_{ab} q^\beta e^{\beta t_1}
\frac{t_2^at_3^b}{a!b!} , & g=0 , \\
\displaystyle - \frac{t_1}{4} + \sum_{4\beta=a+2b}
N^{(1)}_{ab} q^\beta e^{\beta t_1} \frac{t_2^at_3^b}{a!b!} , & g=1 , \\
\displaystyle c_g + \sum_{4\beta=a+2b} N^{(g)}_{ab} q^\beta e^{\beta t_1}
\frac{t_2^at_3^b}{a!b!} , & g>1 .
\end{cases}$$
Here, $t_i$ is the formal variable of degree $2i-2$ dual to $\om^i\in
H^{2i}(\CP^3,\Q)$ and $q$ is the generator of the Novikov ring
$\Lambda\cong\Q\(q\)$ of $\CP^3$. By Proposition (3.7), the coefficient of
$t_1$ in $F_1(\CP^3)$ equals $-c_2(\CP^3)/24$, while the rational numbers
$c_g$ are related to the intersection theory of $\Mbar_g$.

Thus, the coefficient $N^{(g)}_{ab}$ is a rational number which ``counts''
the number of stable maps of degree $\beta$ from a curve of genus $g$ to
$\CP^3$ meeting $a$ generic lines and $b$ generic points.

It is shown by Fulton and Pandharipande \cite{FP} that $N^{(0)}_{ab}$
equals the number of rational space curves of degree $\beta$ which meet $a$
generic lines and $b$ generic points. In particular, they are non-negative
integers. By contrast, the coefficients $N^{(1)}_{ab}$ are neither positive
nor integral: for example, $N^{(1)}_{02}=-1/12$. In \cite{cp3}, we prove
the following result.
\begin{theorem}
The number of elliptic space curves of degree $\beta$ passing through $a$
generic lines and $b$ generic points, where $4\beta=a+2b$, equals
$N^{(1)}_{ab} + (2\beta-1)N^{(0)}_{ab}/12$.
\end{theorem}

By evaluating the equation of Proposition (3.14) on
$\om\boxtimes\om\boxtimes\om\boxtimes\om$, we obtain the following relation
among the elliptic Gromov-Witten for $\CP^3$: if $a\ge2$, then
\begin{multline*}
3 N^{(1)}_{ab} = 4 nN^{(1)}_{a-2,b+1} - \tfrac{1}{4} n^2 N^{(0)}_{ab} +
\tfrac{1}{6} n^3 (n-3) N^{(0)}_{a-2,b+1} \\
\shoveleft{ {} - 2 \sum_{\substack{a-2=a_1+a_2\\b+1=b_1+b_2}}
\textstyle N^{(1)}_{a_1b_1} N^{(0)}_{a_2b_2} n_2^2 (n-3n_1)
\binom{a-2}{a_1} \Bigr\{ n_1 \binom{b}{b_1} + n_2 \binom{b}{b_1-1} \Bigr\} } \\
\shoveleft{ {} + \sum_{\substack{a=a_1+a_2\\b=b_1+b_2}} N^{(1)}_{a_1b_1}
N^{(0)}_{a_2b_2} \textstyle \Bigl\{ n_1n_2 (n+3n_1) \binom{a-2}{a_1} +
n_2^2 (3n_1-n) \binom{a-2}{a_1-1} - 6 n_2^3 \binom{a-2}{a_1-2} \Bigr\}
\binom{b}{b_1} } \\
\shoveleft{{} + \tfrac{1}{12} {\displaystyle
\sum_{\substack{a=a_1+a_2\\b=b_1+b_2}}} N^{(0)}_{a_1b_1} N^{(0)}_{a_2b_2}
n_1 n_2^2 } \\
\textstyle \Bigl\{ n_1^2 (3-n_1) \binom{a-2}{a_1} + n_1n_2(n-3n_1-3)
\binom{a-2}{a_1-1} + n_2^2 (-n_1+n_2-6) \binom{a-2}{a_1-2} \Bigr\}
\binom{b}{b_1} \\
\shoveleft{{} + \tfrac{1}{2} \sum_{\substack{a=a_1+a_2+a_3\\b=b_1+b_2+b_3}}
\textstyle N^{(1)}_{a_1b_1} N^{(0)}_{a_2b_2} N^{(0)}_{a_3b_3} \Bigl\{
2n_1n_2^3n_3(n+3n_1-3n_2) \binom{a-2}{a_2,a_3-2} - 6 n_2^3n_3^3
\binom{a-2}{a_2,a_3} } \\
\textstyle {} + n_2^2n_3^2 (3n_1-n) \Bigl( n_1
\binom{a-2}{a_2-1,a_3-1} + n_2 \binom{a-2}{a_2,a_3-1} + n_3
\binom{a-2}{a_2-1,a_3} \Bigr) \Bigr\} \binom{b}{b_2,b_3} .
\end{multline*}
This relation determines the elliptic coefficient $N^{(1)}_{ab}$ for $a>0$
in terms of $N^{(1)}_{0,\frac{1}{2}a+b}$, the elliptic coefficients of
lower degree, and the rational coefficients. To determine $N^{(1)}_{0,b}$,
we need the relation obtained by evaluating Proposition (3.14) on
$\om^2\boxtimes\om^2\boxtimes\om\boxtimes\om$: if $b\ge2$, then
\begin{multline*}
0 = N^{(1)}_{ab}  + \tfrac{1}{24} n(2n-1) N^{(0)}_{a+2,b-1}
+ \tfrac{1}{48} N^{(0)}_{a+4,b-2} \\
\shoveleft{ {} + \sum_{\substack{a+2=a_1+a_2\\b-1=b_1+b_2}} \textstyle
N^{(1)}_{a_1b_1} N^{(0)}_{a_2b_2} \textstyle \Bigl\{ n_2 \Bigl( n
\binom{a}{a_1} + n_2 \binom{a}{a_1-1} \Bigr) \binom{b-2}{b_1-1} } \\[-10pt]
\shoveright{ \textstyle {} - \frac{1}{6} \Bigl( n_1(6n_1-n_2)
\binom{a}{a_1} + n_2 (16n_1-n_2) \binom{a}{a_1-1} + 6n_2^2 \binom{a}{a_1-2}
\Bigr) \binom{b-2}{b_1} \Bigr\} } \\
\shoveleft{{} - \tfrac{1}{12} \sum_{\substack{a+4=a_1+a_2\\b-2=b_1+b_2}} \textstyle
N^{(1)}_{a_1b_1} N^{(0)}_{a_2b_2} \Bigl( n_1 \binom{a}{a_1} + (2n_1-5n_2)
\binom{a}{a_1-1} + 6n_2 \binom{a}{a_1-2} \Bigr) \binom{b-2}{b_1} } \\
\shoveleft{ {} - \tfrac{1}{48} \sum_{\substack{a+4=a_1+a_2\\b-2=b_1+b_2}}
N^{(0)}_{a_1b_1} N^{(0)}_{a_2b_2} \textstyle \Bigl( n_1^3(n_1-1) \binom{a}{a_1}
+ n_1^2n_2(2n_1-2n_2+1) \binom{a}{a_1-1} } \\
\textstyle {} + n_1n_2^2(2n_1-2n_2+7) \binom{a}{a_1-2}
+ n_2^3(2n_1+5) \binom{a}{a_1-3} + n_2^4 \binom{a}{a_1-4} \Bigr) \binom{b-2}{b_1} \\
\shoveleft{{} - \tfrac{1}{12} \sum_{\substack{a+4=a_1+a_2+a_3\\b-2=b_1+b_2+b_3}}
\textstyle N^{(1)}_{a_1b_1} N^{(0)}_{a_2b_2} N^{(0)}_{a_3b_3}
\textstyle \Bigl\{ 3n_2n_3 \Bigl( n_2^2 \binom{a}{a_2,a_3-2} + n_3^2
\binom{a}{a_2-2,a_3} \Bigr) } \\
\shoveleft{ {} + \textstyle n_1 \Bigl( n_2^3 \binom{a}{a_2,a_3-4}
+ n_2^2(6n_1-n_3) \binom{a}{a_2-1,a_3-3} - 7n_2n_3^2 \binom{a}{a_2-2,a_3-2}
- 5n_3^3 \binom{a}{a_2-3,a_3-1} \Bigr) } \\
\shoveleft{\textstyle {} + \Bigl( n_2^3(n_1-5n_3) \binom{a}{a_2,a_3-3} 
+ n_2^2n_3(5n_1-7n_3) \binom{a}{a_2-1,a_3-2} } \\
\textstyle {}+ n_2n_3^2 (5n_1-n_3) \binom{a}{a_2-2,a_3-1}
+ n_3^3(n_1+n_3) \binom{a}{a_2-3,a_3}\Bigr) \Bigr\} \binom{b-2}{b_2,b_3} .
\end{multline*}
This relation determine the coefficient $N^{(1)}_{0b}$ in terms of elliptic
coefficients of lower order and the rational coefficients, and thus
ultimately in terms of $N^{(0)}_{02}=1$, the number of lines between two
points.

Using these relation, we obtain the results of Table 1. Up to degree $3$,
Theorem A is easily seen to hold, since there are no elliptic space curves
of degrees $1$ and $2$, while all elliptic space curves of degree $3$ lie
in a plane.

It is well-known that there is one quartic elliptic space curve through $8$
general points, while the number of elliptic quartic space curves through
$16$ general lines was calculated by Vainsencher and Avritzer
(\cite{Vainsencher}; see also \cite{Avritzer}, which contains a correction
to \cite{Vainsencher}, bringing it into agreement with our calculation!).

\begin{table} \label{CP3}
\caption{Rational and elliptic Gromov-Witten invariants of $\CP^3$}
$$\begin{tabular}{|R|C|R|D{.}{}{2}|R|} \hline
n & (a,b) & N^{(0)}_{ab} & N^{(1)}_{ab} &
{\scriptstyle N^{(1)}_{ab}+(2n-1)N^{(0)}_{ab}/12} \\ \hline
1 & (0,2) & 1 & -.\frac{1}{12} & 0 \\
& (2,1) & 1 & -.\frac{1}{12} & 0 \\
& (4,0) & 2 & -.\frac{1}{6} & 0 \\[5pt]
2 & (0,4) & 0 & .0 & 0 \\
& (2,3) & 1 & -.\frac{1}{4} & 0 \\
& (4,2) & 4 & -1. & 0 \\
& (6,1) & 18 & -4.\frac{1}{2} & 0 \\
& (8,0) & 92 & -23. & 0 \\[5pt]
3 & (0,6) & 1 & -.\frac{5}{12} & 0 \\
& (2,5) & 5 & -2.\frac{1}{12} & 0 \\
& (4,4) & 30 & -12.\frac{1}{2} & 0 \\
& (6,3) & 190 & -78.\frac{1}{6} & 1 \\
& (8,2) & 1\,312 & -532.\frac{2}{3} & 14 \\
& (10,1) & 9\,864 & -3\,960. & 150 \\
& (12,0) & 80\,160 & -31\,900. & 1\,500 \\[5pt]
4 & (0,8) & 4 & -1.\frac{1}{3} & 1 \\
& (2,7) & 58 & -29.\frac{5}{6} & 4 \\
& (4,6) & 480 & -248. & 32 \\
& (6,5) & 4\,000 & -2\,023.\frac{1}{3} & 310 \\
& (8,4) & 35\,104 & -17\,257.\frac{1}{3} & 3\,220 \\
& (10,3) & 327\,888 & -156\,594. & 34\,674 \\
& (12,2) & 3259\,680 & -1\,515\,824. & 385\,656 \\
& (14,1) & 34\,382\,544 & -15\,620\,216. & 4\,436\,268 \\
& (16,0) & 383\,306\,880 & -170\,763\,640. & 52\,832\,040 \\[5pt]
5 & (0,10) & 105 & -36.\frac{3}{4} & 42 \\
& (2,9) & 1\,265 & -594.\frac{3}{4} & 354 \\
& (4,8) & 13\,354 & -6\,523.\frac{1}{2} & 3\,492 \\
& (6,7) & 139\,098 & -66\,274.\frac{1}{2} & 38\,049 \\
& (8,6) & 1\,492\,616 & -677\,808. & 441\,654 \\
& (10,5) & 16\,744\,080 & -7\,179\,606. & 5\,378\,454 \\
& (12,4) & 197\,240\,400 & -79\,637\,976. & 68\,292\,324 \\
& (14,3) & 2\,440\,235\,712 & -928\,521\,900. & 901\,654\,884 \\
& (16,2) & 31\,658\,432\,256 & -11\,385\,660\,384. & 12\,358\,163\,808 \\
& (18,1) & 429\,750\,191\,232 & -146\,713\,008\,096. & 175\,599\,635\,328 \\
& (20,0) & 6\,089\,786\,376\,960 & -1\,984\,020\,394\,752. &
2\,583\,319\,387\,968 \\ \hline
\end{tabular}$$
\end{table}

\end{document}